\documentclass[11pt]{article}
\usepackage{amsmath}
\usepackage{amssymb}
\usepackage{amsthm}
\usepackage{cases}
\usepackage{url}

\newtheorem{theorem}{Theorem}[section]

\newtheorem{corollary}[theorem]{Corollary}

\newcommand\sa{\smallskipamount}
\newcommand\sPP{\\[\sa]\indent}
\newcommand\al\alpha
\newcommand\be\beta
\newcommand\de\delta
\newcommand\tha\theta
\newcommand\la\lambda
\newcommand\La\Lambda
\newcommand\ga{\gamma}
\newcommand\Ga{\Gamma}
\begin{document}
\title{On the higher-order differential equations for the generalized Laguerre polynomials and Bessel functions}
\author{Clemens Markett}
\date{}
\maketitle

\numberwithin{equation}{section}
\numberwithin{theorem}{section}
\begin{abstract}

In the enduring, fruitful research on spectral differential equations with polynomial eigenfunctions, Koornwinder's generalized Laguerre polynomials are playing a prominent role. Being orthogonal on the positive half-line with respect to the Laguerre weight and an additional point mass $N \ge 0$ at the origin, these polynomials satisfy, for any $\al\in\mathbb{N}_{0}$, a linear differential equation of order $2\al+4$. In the present paper we establish a new elementary representation of the corresponding 'Laguerre-type' differential operator and show its symmetry with respect to the underlying weighted scalar product. Furthermore, we discuss various other representations of the operator, mainly given in factorized form, and show their equivalence. Finally, by applying a limiting process to the Laguerre-type equation, we deduce new elementary representations for the higher-order differential equation satisfied by the Bessel-type functions on the positive half-line.
\end{abstract}
\section{Introduction}
\label{intro}
In 1984, Koornwinder \cite{Ko} introduced the generalized Jacobi polynomials $\{P_{n}^{\al,\be,M,N}(x)\} _{n=0}^{\infty}\,$, $\al,\be >-1$, $M,N \ge 0$, which are  orthogonal with respect to a linear combination of the Jacobi weight function $w_{\al,\be}$ and, in general, two delta “functions” at the endpoints of the interval \linebreak $-1 \le x \le 1$,
\begin{equation}
\begin{aligned}
&w_{\al,\be,M,N}(x)=w_{\al,\be}(x)+M\de(x+1)+N\de(x-1),\\ 
&w_{\al,\be}(x)=h_{\al,\be}^{-1}(1-x)^{\al}(1+x)^{\be},\\
&h_{\al,\be}=\int_{-1}^{1}(1-x)^{\al}(1+x)^{\be}dx=
2^{\al+\be+1}\Gamma(\al+1)\Ga(\be+1)/\Ga(\al+\be+2).
 \label{eq1.1}
\end{aligned}
\end{equation}

In the particular case $M=0$, $N>0$, the so-called Jacobi-type polynomials are given in terms of the classical Jacobi polynomials by

\begin{equation}
P_n^{\al,\be,0,N}(x)=P_n^{\al,\be}(x)+N\,R_n^{\al,\be}(x)\; n\in\mathbb{N}_{0}=\lbrace 0, 1, \cdots \rbrace,
 \label{eq1.2}
\end{equation} 

where, see \cite {Ba0}, \cite {KK2}, \cite{Ko}, \cite {Ma2}, 
\begin{equation}
\begin{aligned}
&P_n^{\al,\be}(x)=\frac{(\al+1)_n}{n!}{}_2F_1\big(-n,n+\al+\be+1;\al+1;\frac{1-x}{2}\big),\;n\in\mathbb{N}_{0},\\
&R_n^{\al,\be}(x)=r_n^{\al,\be}(x-1)P_{n-1}^{\al+2,\be}(x),\; r_n^{\al,\be}=\frac{(\al+\be+2)_n(\al+2)_{n-1}}{2n!\,(\be+1)_{n-1}}, n \ge 1,\;R_0^{\al,\be}(x)=0.
 \label{eq1.3}
\end{aligned}
 \end{equation}	  
The Jacobi-type polynomials, in turn, give rise to the Laguerre-type polynomials $\{L_{n}^{\al,N}(x)\} _{n=0}^{\infty}\,$, $\al >-1,\;N>0$, via the “confluent” limit

\begin{equation}
L_{n}^{\al,N}(x)=\lim_{\be \to \infty}P_n^{\al,\be,0,N}\big(1-\frac{2x}{\be}\big)=L_n^{\al}(x)+N\,T_n^{\al}(x),\; 0 \le x < \infty.
 \label{eq1.4}
\end{equation} 
Here, $L_{n}^{\al}(x)=[(\al +1)_n/n!]\;{}_1F_1(-n;\al+1;x),\;n\in\mathbb{N}_{0}$, denote the Laguerre polynomials, while 
\begin{equation}
T_{n}^{\al}(x)=-t_{n}^{\al}\,x\,L_{n-1}^{\al +2}(x),\;
t_{n}^{\al}=(\al +2)_{n-1}/n!,\; n \ge 1,\;T_0^{\al}(x)=0.
 \label{eq1.5}
\end{equation} 	

In the present paper we are interested in two characteristic features of the Laguerre-type polynomials for any $\al\in\mathbb{N}_{0}$ and $N>0$. Firstly, they are orthogonal with respect to the scalar product induced by the Laguerre weight function $w_{\al}(x)=\al!^{-1} e^{-x}x^{\al}$ plus a delta distribution at the origin, i. e.,
\begin{equation}
(f,g)_{w(\al,N)}= \int_{0}^{\infty}f(x)g(x)w_{\al}(x)dx+N\,f(0)g(0),\;
f,g \in C[0,\infty).
 \label{eq1.6}
\end{equation}  
And secondly, they satisfy a spectral differential equation of order $2\al+4$, which may be stated in operational form by 
  \begin{equation}
  \bigg\lbrace \big\lbrack L_{2,x}^{\al}+n\big\rbrack +  
   \frac{N}{(\al+2)!}\big\lbrack L_{2\al+4,x}^{\al}+(n)_{\al +2}\big\rbrack
    \bigg\rbrace L_n^{\al,N}(x)=0,\;n\in\mathbb{N}_{0},\;0 \le x<\infty.
        \label{eq1.7}
    \end{equation}
While for $N=0$, Eq. (\ref{eq1.7}) reduces to the classical Laguerre equation 

 \begin{equation}
  L_{2,x}^{\al}L_n^{\al}(x)=-n\,L_n^{\al}(x),\;
   L_{2,x}^{\al}=xD_x^2+(\al +1-x)\,D_x,\;n\in\mathbb{N}_{0},\;0 \le x<\infty,
 \label{eq1.8}
    \end{equation}
J. and R. Koekoek \cite [(15)]{KK1} used some deeper special function arguments to determine the higher-order component of the differential operator for sufficiently smooth functions $y(x)$ by 
 \begin{equation}
 \begin{aligned}
 L_{2\alpha+4,x}^{\al}y(x)=&\sum_{i=1}^{2\al+4}d_i^{\al}(x)D_x^i y(x),\\
d_i^{\al}(x)=&\sum_{j=\max(1,\,i-\al-2)}^{\min(i,\al+2)}
(-1)^{i+j+1}\binom{\al+1}{j-1}\binom{\al+2}{i-j}(i+1)_{\al+2-j}\,x^j.
\label{eq1.9}
\end{aligned}
  \end{equation}
Here and in the following, $D_x^i\equiv (D_x)^i$, $i=1,2,\cdots$, denotes an $i$-fold differentiation with respect to $x$. Moreover, we indicate the order of a differential operator as another index.

For the lowest values of the parameter, $\al =0, 1, 2$, the respective fourth-, sixth- and eighth-order Laguerre-type equations have already been established by Littlejohn and Krall, see e.g. \cite{LK}. Shortly after the discovery of the Laguerre-type operator (\ref{eq1.9}), Everitt, Littlejohn, and Wellman \cite{ELW}, \cite{We} found its Lagrange symmetric form, 
 \begin{equation}
 \begin{aligned}
 L_{2\alpha+4,x}^{\al}y(x)=&e^x x^{-\al} \sum_{k=1}^{\al+2}(-1)^{k+1}
  D_x^k \big\lbrack b_k^{\al}(x) D_x^k y(x) \big\rbrack,\\
b_k^{\al}(x)=&\frac{(\al+1)!\;(\al+2)!}{(k-1)!\;k!}
\sum_{j=2k-2}^{k+\al}\frac{e^{-x} x^j}{(j-2k+2)!}.
\label{eq1.10}
\end{aligned}
  \end{equation}

The first major purpose of this paper is to state the operator $L_{2\alpha+4,x}^{\al}$ in another, completely elementary form, see Theorem 2.1 below. This result is then used to show that the components (\ref{eq1.5}) of the Laguerre-type polynomials themselves are eigenfunctions of the second part of equation (\ref{eq1.7}), i. e.
 \begin{equation}
 \big\lbrack L_{2\al+4,x}^{\al}+(n)_{\al +2}\big\rbrack
 T_n^{\al}(x)=0,\;n\in\mathbb{N}_{0}.
 \label{eq1.11}
    \end{equation}
  
 One way to achieve the new representation of the Laguerre-type operator $L_{2\alpha+4,x}^{\al}$ is to apply the limiting process (\ref{eq1.4}) to the higher-order differential equation for the Jacobi-type polynomials (\ref{eq1.2}) which we recently established for all $\al \in \mathbb{N}_0$, $\be >-1, N>0$ \cite{Ma2}. Alternatively, when inserting the definition of the Laguerre-type polynomials into equation (\ref{eq1.7}) for arbitrary  $N>0$ and taking account of the identities (\ref{eq1.8}) and (\ref{eq1.11}), a direct proof of Theorem \ref{thm2.1} depends on showing that 
  
  \begin{equation}
  (\al+2)!\;\lbrack L_{2,x}^{\al}+n\big\rbrack T_n^{\al}(x)+ \lbrack  
   L_{2\al+4,x}^{\al}+(n)_{\al +2}\big\rbrack L_n^{\al}(x)=0,\;n\in\mathbb{N}_{0}.
  \label{eq1.12}
    \end{equation}
This, however, is an immediate consequence of the two identities stated in Theorem \ref{thm2.3}. Finally, once knowing the new representation (2.5), it is not hard to verify its equivalence to Koekoek's representation (\ref{eq1.9}), see Corollary 2.4.

Another central goal of this paper is pursued in Section 3. Here we deal with four additional representations of $L_{2\alpha+4,x}^{\al}$, each being given as a product of $\al+2$ second-order differential expressions. One of them is a new, completely non-commutative factorization which enables us to set up a recurrence relation for the Laguerre-type operator with respect to the parameter $\al$, see Theorem 3.1, Eq. (\ref{eq3.5}). The other three factorizations are essentially due to Bavinck \cite [Secs. 2.1 and 2.2]{Ba0} and Dur\'{a}n \cite{Du}. It is remarkable to observe how these representations are related to each other. In particular, we verify the characteristic property (\ref{eq1.11}) for some of them, see Corollaries 3.2 and 3.3 (ii). 

In Section 4 we show that the combined Laguerre-type differential operator

 \begin{equation}
  L_{2\al+4,x}^{\al,N}y(x)= L_{2,x}^{\al}y(x) +
  N\;(\al +2)!^{-1} L_{2\al+4,x}^{\al}y(x), \;y\in C^{(2\al +4)}[0,\infty),
  \label{eq1.13}
    \end{equation}
is symmetric with respect to the scalar product (\ref{eq1.6}). Consequently, its eigenfunctions form an orthogonal system in the corresponding function space and thus recover the known orthogonality of the Laguerre-type polynomials. By the way, one obtains the same result when utilizing the Lagrange symmetric form (\ref{eq1.10}) of the Laguerre-type operator. It would be worthwhile to compare this former representation with our new one within a deeper spectral theoretical context. 

In the last Section 5, we establish new elementary representations of the higher-order Bessel-type differential equation as well. This equation has been introduced in the joint paper with W. N. Everitt \cite{EM}. Its eigenfunctions, the so-called Bessel-type functions $\{J_{\la}^{\al,M}(r)\} _{\la \ge 0}$, $\al \in \mathbb{N}_0$, $M \ge 0$, $0 \le r<\infty$, satisfy an orthogonality relation in distributional sense and thus give rise to a generalized Hankel transform and its inverse, cf. \cite{EKLM}. In the present paper we basically use the fact that the Bessel-type functions are determined via another confluent limiting process, now originating in the Laguerre-type polynomials. What is left is to examine how the properties of the Laguerre-type equation carry over to the Bessel-type case. 

From a more general point of view, there is a vast literature on orthogonal polynomial systems satisfying a finite order differential equation of a certain spectral type. For a characterization of these so-called Bochner-Krall polynomials and some major contributions in the field see Kwon, Littlejohn, and Yoon \cite{KLY}. In this context, our new results on the Laguerre-type equation may serve as a prototype to investigate other basic phenomena in technically more involved situations. Only recently, we established an elementary representation of the differential equation of order $2\al+2\be+6$ for the generalized Jacobi polynomials $\{P_{n}^{\al,\be,M,N}(x)\} _{n=0}^{\infty}\,$, $\al, \be \in \mathbb{N}_0$, $M, N>0$ \cite{Ma3}. Moreover, one may think of possible extensions to differential equations of Sobolev-type or even to higher-order difference equations satisfied by discrete orthogonal polynomials.

   \section{New elementary representation of the Laguerre-type equation}
 \label{sec:2}
 
In the paper \cite{Ma2}, we stated the higher-order differential equation for the Jacobi-type polynomials $P_n^{\al,\be,0,N}(x)$, $\al \in \mathbb{N}_0, \be >-1$, $N>0$, in the form
  \begin{equation}
 \bigg\lbrace \big\lbrack L_{2,x}^{\al,\be}-\La_{2,n}^{\al,\be}\big\rbrack +  
   \frac{N}{b_{\al,\be}}\big\lbrack L_{2\al+4,x}^{\al,\be}
         -\La_{2\al+4,n}^{\al,\be}\big\rbrack \bigg\rbrace P_n^{\al,\be,0,N}(x)=0,
  \; n\in \mathbb{N}_0,
        \label{eq2.1}
    \end{equation}
where, for sufficiently smooth functions $y(x)$, the two differential expressions are given by 
 \begin{equation}
  \begin{aligned}
   L_{2,x}^{\al,\be}y(x)=&\big\lbrace(x^2-1)D_x^2+\left[\al-\be+(\al+\be+2)x\right]
   D_x \big\rbrace y(x)\\
   =&(x-1)^{-\al}(x+1)^{-\be}D_x \big\lbrack (x-1)^{\al+1}(x+1)^{\be +1} y(x)\big\rbrack,
    \end{aligned}
     \label{eq2.2}
      \end{equation} 
 \begin{equation} 
  L_{2\al +4,x}^{\al ,\be}y(x) =\frac{x-1}{(x+1)^{\be}}D_x^{\al+2}\big\lbrace (x+1)^{\al +\be+2}D_x^{\al+2} \big\lbrack (x-1)^{\al+1} y(x)\big\rbrack \big\rbrace.
  \label{eq2.3}
 \end{equation}     
The two components in Eq. (\ref{eq2.1}) are linked by the constant $b_{\al,\be}=(\al +2)!\;(\be+1)_{\al+1}$, and the eigenvalue parameters are monic polynomials in $n$  of order 2 and $2\al+4$, respectively, namely 
 \begin{equation}
    \La_{2,n}^{\al,\be}=n(n+\al+\be+1),\;
    \La_{2\al+4,n}^{\al,\be}=(n)_{\al+2}\,(n+\be)_{\al+2}\,.
     \label{eq2.4}
    \end{equation}

Applying now the limit relation (\ref{eq1.4}) to the Jacobi-type equation (2.1--4), we arrive at the following fundamental result. 
 \begin{theorem}
 \label{thm2.1}
 Let $\al \in \mathbb{N}_0$, $N>0$. Koornwinder's Laguerre-type polynomials $L_n^{\al,N}(x)$, $n \in \mathbb{N}_0$, satisfy the eigenvalue equation ((\ref{eq1.7}), where, for sufficiently smooth functions $y(x)$,
  \begin{equation}
   \begin{aligned}
   L_{2,x}^{\al}y(x)&=e^x x^{-\al}D_x\big\lbrace e^{-x} x^{\al+1}D_x y(x)\big\rbrace,\\
   L_{2\al+4,x}^{\al}y(x)&=(-1)^{\al+1} e^x x\,D_x^{\al+2}\big\lbrace e^{-x}D_x^{\al+2}[x^{\al+1}y(x)]\big\rbrace .
   \end{aligned}
  \label{eq2.5}
  \end{equation}
   \end{theorem}
 \begin{Proof}
 For convenience, let us choose a new independent variable $\xi$ in Eq. (\ref{eq2.1}) and divide the equation by $-\be$. Hence the Jacobi-type polynomials satisfy 
  \begin{equation}
  \bigg\lbrace \frac{1}{\be}\big\lbrack -L_{2,\xi}^{\al,\be}+\La_{2,n}^{\al,\be}
  \big\rbrack +\frac{N}{\be\,b_{\al ,\be}}\big\lbrack -L_{2\al+4,\xi}
  ^{\al,\be}+\La_{2\al+4,n}^{\al,\be}\big\rbrack\bigg\rbrace P_n^{\al,\be,0,N}(\xi)=0.
  \label{eq2.6}
  \end{equation}   
 In the limit $\be\to\infty$, the eigenvalue parameters tend to
 \begin{equation*}
 \frac{\La_{2,n}^{\al,\be}}{\be}=\frac{n(n+\al+\be+1)}{\be}\to n,\;
 \frac{ \La_{2\al+4,n}^{\al,\be}}{\be\,b_{\al,\be}}=
 \frac{(n)_{\al +2}\,(n+\be)_{\al +2}}{\be(\al+2)!\,(\be+1)_{\al +1}}\to
 \frac{(n)_{\al +2}}{(\al +2)!}.
 \end{equation*}
 What the differential expressions in Eq. (2.6) are concerned, we substitute $\xi=1-2x/\be$ and thus replace $D_\xi$ formally by $-(\be /2)D_x$. Hence,
 \begin{equation*}
 \begin{aligned}
  - \frac{1}{\be}\,&L_{2,\xi}^{\al,\be}P_n^{\al,\be ,0,N}(\xi)\\
  &=\bigg\lbrace \bigg\lbrack \frac{4x}{\be}-\big(\frac{2x}{\be}\big)^2 \bigg\rbrack\frac{\be}{4}D_x^2+\bigg\lbrack 2\al+2-(\al+\be+2)\frac{2x}{\be}
  \bigg\rbrack \frac{1}{2}D_x\bigg\rbrace P_n^{\al ,\be ,0,N}\bigg(1-\frac{2x}
  {\be}\bigg)\\
  &\xrightarrow[\be \to \infty]{} \big\lbrack xD_x^2+(\al+1-x)D_x\big\rbrack L_n^{\al,N}(x)
  =e^x x^{-\al}D_x\big\lbrack e^{-x} x^{\al+1}D_x L_n^{\al,N}(x)\big\rbrack
    \end{aligned}
  \end{equation*}   
and
\begin{equation*}
\begin{aligned}
  -&\frac{1}{\be\,b_{\al,\be}}\,L_{2\al+4,\xi}^{\al,\be}P_n^{\al,\be,0,N}(\xi)=
  \frac{2x}{\be^2(\al+2)!\,(\be+1)_{\al+1}}\\
  &\cdot \bigg(2-\frac{2x}{\be}\bigg)^{-\be} \bigg(\frac{\be}{2}\bigg)^{2\al+4}
  D_x^{\al+2}\bigg\lbrace \bigg(2-\frac{2x}{\be}\bigg)^{\al +\be +2}D_x^{\al+2}
  \bigg\lbrack\bigg(-\frac{2x}{\be}\bigg)^{\al+1}P_n^{\al,\be,0,N}
  \bigg(1-\frac{2x}{\be}\bigg)\bigg\rbrack\bigg\rbrace\\
  &\xrightarrow[\be \to \infty]{} (-1)^{\al+1}\frac{e^x x}{(\al+ 2)!}D_x^{\al+2}\big\lbrace e^{-x}D_x^{\al+2}\big\lbrack x^{\al+1}L_n^{\al ,N}(x)\big\rbrack\big\rbrace.
    \end{aligned}
    \end{equation*}
All limits can clearly be made rigorous, so that Eq. (2.6) turns into the Laguerre-type equation (\ref{eq1.7}) with $L_{2,x}^\al$ and $L_{2\al+4,x}^\al$ as given in item (\ref{eq2.5}). \hfill $\square$    
\end{Proof}

\begin{corollary}
\label{cor2.2}
For $\al \in \mathbb{N}_0$ let the operator $L_{2\al+4,x}^\al$ be defined as in item (\ref{eq2.5}). Then for all $n \in \mathbb{N}$, the components $T_n^\al(x)$ of the Laguerre-type polynomials (\ref{eq1.4}) satisfy the eigenvalue equation (\ref{eq1.11}).
\end{corollary}

\begin{Proof}
We iteratively use the two differentiation formulas for the Laguerre polynomials  (cf. (2.10))
\begin{equation}
\begin{aligned}
D_x\big\lbrack e^{-x}L_n^\ga(x)\big\rbrack &=-e^{-x}L_n^{\ga+1}(x),\;\ga >-1,\\
D_x\big\lbrack x^\ga L_n^\ga(x)\big\rbrack &=(n+\ga)x^{\ga-1}L_n^{\ga-1}(x),
\;\ga >0.
 \end{aligned}
 \label{eq2.7}
 \end{equation}
Then, by definition (\ref{eq1.5}) of $T_n^\al(x)$, it follows that
\begin{equation*}
\begin{aligned}
L_{2\al+4,x}^\al T_n^\al(x) &=(-1)^\al t_n^\al\,e^x x D_x^{\al+2}
\big\lbrace e^{-x}D_x^{\al+2}\big\lbrack x^{\al+2}L_{n-1}^{\al+2}(x)
\big\rbrack\big\rbrace\\
&=(-1)^\al t_n^\al\,e^x x D_x^{\al+2}\big\lbrace e^{-x}(n)_{\al+2}\,
L_{n-1}^0(x)\big\rbrack\big\rbrace\\
&=t_n^\al\,x\,(n)_{\al+2}\,L_{n-1}^{\al+2}(x)=-(n)_{\al+2}\,T_n^\al(x). \hspace{2cm} \square 
 \end{aligned}
 \end{equation*}
\end{Proof}

As carried out in the Introduction, we can prove Theorem 2.1 directly by combining the following two identities to verify the required identity (\ref{eq1.12}).

\begin{theorem}
 \label{thm2.3}
For $\al \in \mathbb{N}_0$ let the operators $L_{2,x}^\al$ and $L_{2\al+4,x}^\al$ be defined in item (2.5). For all $n \in \mathbb{N}$, the two components of the Laguerre-type polynomials (1.4) satisfy  

  \begin{equation}
   (\al+2)!\;\big\lbrack L_{2,x}^{\al}+n\big\rbrack T_n^\al(x)=
        -(n+1)_\al (\al +1)(\al +2)\,L_{n-1}^{\al+2}(x),
         \label{eq2.8}
  \end{equation}
  \begin{equation}
    \big\lbrack L_{2\al+4,x}^{\al}+(n)_{\al +2}\big\rbrack
    L_n^\al(x)=(n+1)_\al (\al +1)(\al +2)\,L_{n-1}^{\al+2}(x).
   \label{eq2.9}
    \end{equation}
        
\noindent If $n=0$, the expressions on the left-hand sides of the identities (2.8--9) clearly vanish.
   \end{theorem}

\begin{Proof}
Identity (2.8) is a direct consequence of the Laguerre equation (\ref{eq1.8}). In fact,
\begin{equation*}
   \big\lbrack L_{2,x}^\al-(\al+1)/x+n\big\rbrack \,\big\lbrack x\, L_{n-1}^{\al+2}(x)\big\rbrack=x\, \big\lbrack L_{2,x}^{\al+2}+n-1\big\rbrack L_{n-1}^{\al+2}(x)=0
     \end{equation*}
and hence, by definition (\ref{eq1.5}),
 \begin{equation*}
   (\al+2)!\;\big\lbrack L_{2,x}^{\al}+n\big\rbrack T_n^\al(x)=
   - (\al+2)!\;t_n^\al (\al +1)L_{n-1}^{\al+2}(x)=
        -(n+1)_\al (\al +1)(\al +2)\,L_{n-1}^{\al+2}(x).
   \end{equation*}
 Concerning Eq. (\ref{eq2.9}) we need, besides the differentiation formulas (\ref{eq2.7}), the following well-known identities for the Laguerre polynomials \cite[10.12 (15),(23),(24)]{HTF2}
 \begin{equation}
 \begin{aligned}
D_xL_n^\ga(x) =&-L_{n-1}^{\ga +1}(x),\;\gamma>-1,\\
L_n^\ga(x)=&L_n^{\ga+1}(x)-L_{n-1}^{\ga+1}(x),\;\ga >-1,\\
x\,L_n^\ga(x)=&(n+\ga)L_n^{\ga-1}(x)-(n+1)L_{n+1}^{\ga-1}(x),\;\ga >0.
  \end{aligned}
  \label{eq2.10}
  \end{equation}
 Altogether, we obtain
  \begin{equation*}
  \begin{aligned}
L_{2\al+4,x}^\al L_n^\al(x)=&(-1)^{\al+1} e^x x\,D_x^{\al+2}
\big\lbrace e^{-x}D_x^{\al+2}\big\lbrack x^{\al+1}(L_n^{\al+1}(x)
-L_{n-1}^{\al+1}(x))\big\rbrack\big\rbrace\\
=&(-1)^{\al} e^x x\,D_x^{\al+2}
\big\lbrace e^{-x}\big\lbrack (n+1)_{\al +1} L_{n-1}^1(x)-
(n)_{\al +1} L_{n-2}^1(x)\big\rbrack\big\rbrace\\
=&(n+1)_{\al+1}\,x\,L_{n-1}^{\al+3}(x)-
(n)_{\al +1}\,x\,L_{n-2}^{\al+3}(x)\\
=&(n+1)_{\al+1}\big\lbrack (\al+2)L_{n-1}^{\al+2}(x)-n\,L_n^{\al+1}(x)\big\rbrack-\\
&\quad -(n)_{\al +1}
\big\lbrack (\al +2)L_{n-1}^{\al+2}(x)- (n+\al +1)\,L_{n-1}^{\al+1}(x)\big\rbrack\\
=&\big\lbrack (n+1)_{\al+1}-(n)_{\al+1}\big\rbrack(\al+2)L_{n-1}^{\al+2}(x)-
(n)_{\al+2}\big\lbrack L_n^{\al+1}(x)-L_{n-1}^{\al+1}(x)\big\rbrack \\
=&(n+1)_\al(\al+1)(\al+2)L_{n-1}^{\al+2}(x)-
(n)_{\al+2} L_n^\al(x).
   \end{aligned}
    \end{equation*}
\end{Proof}
    \hfill $\square$
    
We close this section by showing the following relationship.

\begin{corollary}
\label{cor2.4}
The representations (\ref{eq1.9}) and (\ref{eq2.5}) of the differential operator $L_{2\al+4,x}^\al$ are equivalent. 
\end{corollary}

\begin{Proof}
Carrying out the two derivatives in item (\ref{eq2.5}), we find that
\begin{equation*}
\begin{aligned}
(-1)^{\al+1}&e^x x\, D_x^{\al +2}\big\lbrace e^{-x}D_x^{\al+2}\big\lbrack x^{\al+1}y(x)\big\rbrack\big\rbrace\\
=&(-1)^{\al+1}\,x\sum_{k=0}^{\al+2}\binom{\al+2}{k}(-1)^{\al+2-k}D_x^{k+\al+2}
\big\lbrack x^{\al+1}y(x)\big\rbrack\\
=&x\sum_{k=0}^{\al+2}\binom{\al+2}{k}(-1)^{k+1}\sum_{i=k+1}^{k+\al+2}
\binom{k+\al+2}{i}\frac{(\al+1)!}{(i-k-1)!}x^{i-k-1}D_x^i y(x)\\
=&\sum_{i=1}^{2\al+4}d_i^\al(x)D_x^i y(x),
    \end{aligned}
  \end{equation*}
where
\begin{equation*}
\begin{aligned}
d_i^\al(x)=&\sum_{k=\max (0,i-\al-2)}^{\min(i-1,\al+2)}\binom{\al+2}{k}
(-1)^{k+1}\binom{k+\al+2}{i}\frac{(\al+1)!}{(i-k-1)!}x^{i-k}\\
=&\sum_{j=\max(1,i-\al-2)}^{\min(i,\al+2)}(-1)^{i+j+1}\binom{\al+2}{i-j}
\binom{\al+1}{j-1}(i+1)_{\al+2-j}\,x^j.\hspace{2cm} \square 
    \end{aligned}
  \end{equation*}
 \end{Proof}

   \section{Factorizations of the Laguerre-type differential operator}
 \label{sec:3}

Recently we investigated two factorizations of the Jacobi-type differential operator into products of second-order differential expressions based on the classical Jacobi differential operator (\ref{eq2.2}), see \cite [Sec. 4]{Ma2}. For sufficiently smooth functions $y(x) \equiv (x-1)u(x)$, both factorizations are given in two equivalent versions by
\begin{equation}
   \begin{aligned}
   L_{2\al+4,x}^{\al,\be}y(x)=&\prod_{j=0}^{\al+1}\bigg\lbrace L_{2,x}^{\al,\be}-\frac{2(\al+1)}{x-1}+j(\al+\be+1-j)\bigg\rbrace y(x),\\
  L_{2\al+4,x}^{\al,\be}[(x-1)u(x)]  =&(x-1)\prod_{j=0}^{\al+1}\big\lbrace L_{2,x}^{\al+2,\be}+(j+1)(\al+\be+2-j)\big\rbrace u(x)
  \label{eq3.1}
    \end{aligned}
     \end{equation}
 and

 \begin{equation}
    \begin{aligned}
    L_{2\al+4,x}^{\al,\be}y(x)=&\prod_{j=0}^{\al+1}\bigg\lbrace L_{2,x}^{2j-1,\be}-\frac{4j}{x-1}+j(j+\be)\bigg\rbrace y(x),\\
   L_{2\al+4,x}^{\al,\be}[(x-1)u(x)]  =&(x-1)\prod_{j=0}^{\al+1}\big\lbrace L_{2,x}^{2j+1,\be}+(j+1)(j+\be+1)\big\rbrace u(x).
   \label{eq3.2}
     \end{aligned}
      \end{equation}
 In the latter identities, the non-commutative, multiple product $\prod_{j=0}^{\al+1}$ is understood in the sense that the second-order differential expressions under the product sign are successively applied to the respective functions on their right-hand side in order from $j=0$  to $j=\al+1$. 
 Analogously, there are two factorizations of the higher-order Laguerre-type operator.

\begin{theorem}
\label{thm3.1}
(i) For any $\al \in \mathbb{N}_0$, the Laguerre-type differential operator $L_{2\al +4,x}^\al$ of order $2\al +4$ can be factorized in two different ways into products of $\al +2$ factors involving the Laguerre differential operator $L_{2,x}^{\ga}=xD_x^2+(\ga +1-x)\,D_x$, $\ga \in \{-1\} \cup \mathbb{N}_0$. In fact, for $y(x) \equiv xu(x) \in C^{(2\al +4)}[0,\infty)$, there hold 
 \begin{equation}
 \begin{aligned}
 L_{2\al +4,x}^\al y(x)=& (-1)^{\al +1}\prod_{j=0}^{\al+1}\bigg\lbrace L_{2,x}^\al -\frac{\al+1}{x}-j\bigg\rbrace y(x),\\
 L_{2\al +4,x}^\al [xu(x)]=&(-1)^{\al +1}x\prod_{j=0}^{\al+1}\big\lbrace L_{2,x}^{\al +2} -j-1\big\rbrace u(x)
    \end{aligned}
     \label{eq3.3}   
   \end{equation}
 and
  \begin{equation}
  \begin{aligned}
  L_{2\al +4,x}^\al y(x)=& (-1)^{\al +1}\prod_{j=0}^{\al+1}\bigg\lbrace L_{2,x}^{2j-1} -\frac{2j}{x}-j\bigg\rbrace y(x),\\
  L_{2\al +4,x}^\al [xu(x)]=&(-1)^{\al +1}x\prod_{j=0}^{\al+1}\big\lbrace L_{2,x}^{2j+1} -j-1\big\rbrace u(x).
     \end{aligned}
      \label{eq3.4}   
    \end{equation}
 (ii) The representations (\ref{eq3.3}) and (\ref{eq3.4}) are equivalent.\\
 (iii) The Laguerre-type differential operator satisfies the recurrence relation
 \begin{equation}
    L_{2\al +4,x}^\al y(x)=-\big\lbrace L_{2,x}^{2\al +1} -(2\al+2)/x-\al -1\big\rbrace L_{2\al +2,x}^{\al-1} y(x),\; \al \in \mathbb{N}_0.
  \label{eq3.5}   
 \end{equation}
    \end{theorem}
 
 \begin{Proof}
  (i) The first identity in item (\ref{eq3.3}) is given already in \cite[(2.1)]{Ba0}. Similarly as in the proof of Theorem 2.1, the representations (3.3--4) follow by applying the limiting process (\ref{eq1.4}) to the factorizations of the Jacobi-type differential operator, (\ref{eq3.1}) and (\ref{eq3.2}), respectively. Concerning item (\ref{eq3.4}), for instance, we use the fact that for all $j\in\mathbb{N}_0$,  
 \begin{equation*}
  \lim_{\be \to \infty} \frac{1}{\be}\bigg\lbrace L_{2,\xi}^{2j-1,\be} -\frac{4j}{\xi-1}+j(j+\be)\bigg\rbrace y(\xi)\,\big\vert_{\xi=1-2x/\be}
  =\bigg\lbrace-L_{2,x}^{2j-1} +\frac{2j}{x}+j\bigg\rbrace \lim_{\be \to \infty} y\bigg(1-\frac{2x}{\be}\bigg).
   \end{equation*}
 
 (ii) In order to verify the equivalence of the two representations, it suffices to compare the second lines of (3.3--4). To this end we are going to show by induction with respect to $\al\in\mathbb{N}_0$,  that for any sufficiently smooth function $u(x)$,
 \begin{equation}
  \prod_{j=1}^{\al+1}\big\lbrace L_{2,x}^{\al +1} -j\big\rbrace u(x)=
   \prod_{j=1}^{\al+1}\big\lbrace L_{2,x}^{2j-1} -j\big\rbrace u(x).
         \label{eq3.6}
       \end{equation}
 While identity (\ref{eq3.6}) is trivial for $\al=0$, we suppose that it holds true for some  $\al\in\mathbb{N}_0$. Concerning the step from $\al$ to $\al+1$, we multiply both sides of (\ref{eq3.6}) from the left by the differential expression $L_{2,x}^{2\al+3}-\al-2$. Hence the new product on the right-hand side runs till $j=\al+2$ as required. On the left-hand side, we use the straightforward commutation relation 
 \begin{equation}
   \big\lbrack L_{2,x}^{\al+j+2} -j-1\big\rbrack\, 
   \big\lbrack L_{2,x}^{\al+1} -j\big\rbrack u(x)=
    \big\lbrack L_{2,x}^{\al+2} -j-1\big\rbrack\, 
      \big\lbrack L_{2,x}^{\al+j+1} -j\big\rbrack u(x),\;j \in \mathbb{N},
          \label{eq3.7}
        \end{equation}
 successively for $j=\al+1, \al, \cdots, 1$, to end up with the product $\prod_{j=1}^{\al+2}\big\lbrace L_{2,x}^{\al+2} -j\big\rbrace u(x)$.\\
  
(iii) The recurrence relation (3.5) is a direct consequence of the first identity in item (\ref{eq3.4}). Notice that for $\al =0$, the identity is in accordance with the representation (\ref{eq1.9}), since
\begin{equation*}
\begin{aligned}
   L_{4,x}^0 y(x)=&-\big\lbrack xD_x^2+(2-x)D_x -2/x-1\big\rbrack\,x\, 
   \big\lbrack D_x^2-D_x\big\rbrack y(x)\\
   =&-x\big\lbrack xD_x^2+(4-x)D_x-2\big\rbrack\, 
      \big\lbrack D_x^2-D_x\big\rbrack y(x)\\
    =&-x\big\lbrack xD_x^4+(4-2x)D_x^3-(6-x)D_x^2+2D_x\big\rbrack y(x).\hspace*{2cm}\square
      \end{aligned}
        \end{equation*} 
    \end{Proof}   

\begin{corollary}
\label{cor3.2}
For $\al \in \mathbb{N}_0$ let the Laguerre-type operator $L_{2\al+4,x}^\al$ be defined via the second product in item (\ref{eq3.3}). Then its eigenfunctions are given by the components $T_n^\al(x)$, $n \in \mathbb{N}$, of the Laguerre-type polynomials (\ref{eq1.4}).
\end{corollary}

\begin{Proof}
See \cite{Ba0}. In view of the Laguerre equation (\ref{eq1.8}) with parameter $\al+2$, it follows that
\begin{equation*}
\begin{aligned}
L_{2\al+4,x}^\al T_n^\al(x)=&(-1)^\al t_n^\al\,x\,\prod_{j=0}^{\al+1} \big\lbrace L_{2,x}^{\al +2} -j-1\big\rbrace L_{n-1}^{\al+2}(x)\\
=&(-1)^\al t_n^\al\,x\,\prod_{j=0}^{\al+1} \big\lbrace -n+1-j-1\big\rbrace \cdot L_{n-1}^{\al+2}(x)\\
=&t_n^\al\,x\,(n)_{\al+2}\,L_{n-1}^{\al+2}(x)=-(n)_{\al+2}\,T_n^\al(x).
\hspace*{2cm}\square
 \end{aligned}
 \end{equation*} 
   \end{Proof}

Recently, A. J. Dur\'{a}n \cite{Du} developed a concept based on so-called $\mathcal{D}$-operators to construct orthogonal polynomial systems that are eigenfunctions of higher-order differential and difference equations. A prominent example deals with the Laguerre-type case. Adjusting the notations used in \cite[Sec. 8.1]{Du} to ours, the orthogonal polynomials are in accordance with the Laguerre-type polynomials (1.4--5). After writing the resulting differential equation in the form (\ref{eq1.7}) and recalling that $L_{2,x}^{-1}=x[D_x^2-D_x]$, we identify the differential operator as the product
\begin{equation}
 L_{2\al +4,x}^\al y(x)=(-1)^{\al +1}L_{2,x}^{-1}\prod_{j=1}^{\al+1}\big\lbrace L_{2,x}^{\al +1} -j\big\rbrace y(x),\; \al \in \mathbb{N}_0.
 \label{eq3.8}   
   \end{equation}

\begin{corollary}
\label{cor3.3}
(i) The representations (\ref{eq3.4}) and (\ref{eq3.8}) are equivalent.\\
(ii) For $\al \in \mathbb{N}_0$, the operator (\ref{eq3.8}) satisfies the eigenvalue equation (\ref{eq1.11}) with eigenfunctions $T_n^\al(x)$, $n\in \mathbb{N}$.
\end{corollary}
\begin{Proof}
(i) In view of the representation (3.4) and the straightforward commutation relation
\begin{equation}
\begin{aligned}
  \bigg\lbrack L_{2,x}^{2j-1}-\frac{2j}{x} -j\bigg\rbrack\,L_{2,x}^{-1}y(x)&=
    x \big\lbrack L_{2,x}^{2j+1} -j-1\big\rbrack\,\big\lbrack D_x^2-D_x\big\rbrack y(x) \\
    &=L_{2,x}^{-1}\big\lbrack L_{2,x}^{2j-1}-j\big\rbrack y(x),\;j \in \mathbb{N},
         \end{aligned}
          \label{eq3.9}
        \end{equation}
  we obtain
     \begin{equation*}
   \begin{aligned} 
   L_{2\al +4,x}^\al y(x)=& (-1)^{\al +1}\prod_{j=1}^{\al+1}\bigg\lbrace L_{2,x}^{2j-1} -\frac{2j}{x}-j\bigg\rbrace L_{2,x}^{-1}\, y(x)\\
      =&(-1)^{\al +1}\,L_{2,x}^{-1}\,\prod_{j=1}^{\al+1}\bigg\lbrace L_{2,x}^{2j-1}-j\bigg\rbrace\,y(x).
    \end{aligned}
    \end{equation*}
 Replacing the last product by means of identity (\ref{eq3.6}) then yields the representation (\ref{eq3.8}).

(ii) Here we successively employ the properties (\ref{eq2.10}) of the Laguerre polynomials and the Laguerre equation (\ref{eq1.8}) with parameter $\al+1$, whence
  \begin{equation*}
    \begin{aligned}
  L_{2\al +4,x}^\al T_n^\al(x)=&(-1)^{\al}t_n^\al\,L_{2,x}^{-1}  \prod_{j=0}^{\al+1} \bigg\lbrace L_{2,x}^{\al +1}-j\bigg\rbrace x\,L_{n-1}^{\al +2}(x) \\
  =&(-1)^{\al}t_n^\al\,L_{2,x}^{-1}  \prod_{j=0}^{\al+1} \bigg\lbrace L_{2,x}^{\al +1}-j\bigg\rbrace \big\lbrack (n+\al+1)L_{n-1}^{\al+1}(x)-n\,L_n^{\al+1}(x)\big\rbrack \\
  =&-t_n^\al\,L_{2,x}^{-1}\bigg\lbrace (n+\al+1)\prod_{j=1}^{\al+1}(n-1+j) \cdot L_{n-1}^{\al+1}(x)-n\prod_{j=1}^{\al+1}(n+j)\cdot L_n^{\al+1}(x)\bigg\rbrace\\       
 =&t_n^\al \,x\,\big\lbrack D_x^2-D_x\big\rbrack (n)_{\al +2}L_n^\al(x)\\
 =&t_n^\al(n)_{\al +2}\,x\,\big\lbrack L_{n-2}^{\al +2}(x)+L_{n-1}^{\al+1}(x)
     \big\rbrack\\
  =&t_n^\al(n)_{\al +2}\,x\,L_{n-1}^{\al +2}(x)=-(n)_{\al +2}T_n^\al(x).\hspace{2cm}\square
      \end{aligned}
      \end{equation*}
    \end{Proof}
 
  Finally, we note that Bavinck \cite[(2.2)]{Ba0} determined another factorization of the Laguerre-type operator which is given, in our notation, by  
 \begin{equation}
  L_{2\al +4,x}^\al y(x)=(-1)^{\al +1}x\,\big\lbrack L_{2,x}^{2\al+3}-\al -2 \big\rbrack \big\lbrack D_x^2-D_x\big\rbrack\,\prod_{j=1}^{\al}\big\lbrack L_{2,x}^{\al} -j\big\rbrack y(x).
  \label{eq3.10}   
    \end{equation}
 
 \begin{corollary}
 \label{cor3.4}
 The representations (3.8) and (3.10) are equivalent.
 \end{corollary}
 \begin{Proof}
 Proceeding from the representation (\ref{eq3.10}), we first apply the commutation relation (\ref{eq3.9}) for $j=\al+1$ to and then, successively for $j=\al,\cdots,1$, the identity (\ref{eq3.7}) with $\al$ replaced by $\al-1$. This yields formula (3.8),
  \begin{equation*}
  \begin{aligned}
   L_{2\al +4,x}^\al y(x)=&(-1)^{\al +1}L_{2,x}^{-1}\,\big\lbrack L_{2,x}^{2\al+1}-\al-1 \big\rbrack\,\prod_{j=1}^{\al}\big\lbrack L_{2,x}^{\al} -j\big\rbrack y(x)\\
   =&(-1)^{\al +1}L_{2,x}^{-1}\,\prod_{j=1}^{\al+1}\big\lbrack L_{2,x}^{\al+1} -j\big\rbrack y(x).\hspace{2cm}\square
  \end{aligned}
   \end{equation*}
       \end{Proof}
    
 \section{Symmetry of the Laguerre-type differential operator and orthogonality of the corresponding eigenfunctions}
  \label{sec:4}    
 
In item (1.13) we introduced the operator $L_{2\al+4,x}^{\al,N}$ as the linear combination (\ref{eq1.13}) of the two differential expressions (\ref{eq2.5}). Hence, the differential equation (\ref{eq1.7}) for the Laguerre-type polynomials $y_n(x)=L_n^{\al,N}(x)$ takes the form 
\begin{equation}
 \bigg\lbrace L_{2\al+4,x}^{\al,N}+\La_{2\al+4,n}^{\al,N}\bigg\rbrace y_n(x)=0 ,\; \La_{2\al+4,n}^{\al,N} =n+\frac{N}{(\al+2)!}(n)_{\al+2}.
 \label{eq4.1} 
  \end{equation}
  
\begin{theorem}
\label{thm4.1}
For $\al \in \mathbb{N}_0$ and $N > 0$, the Laguerre-type differential operator $L_{2\al+4,x}^{\al,N}$ is symmetric with respect to the scalar product (\ref{eq1.6}), i.e.
\begin{equation}
\big(L_{2\al+4,x}^{\al,N}f,g\big)_{w(\al,N)}=\big(f,L_{2\al+4,x}^{\al,N}g \big)_{w(\al ,N)},\;f,g \in C^{(2\al+4)}[0,\infty).
 \label{eq4.2} 
  \end{equation}
  \end{theorem}
\begin{Proof}
The left-hand side of identity (\ref{eq4.2}) is given by
  \begin{equation}
  \begin{aligned}     
  \big(L&_{2\al+4,x}^{\al,N}f,g\big)_{w(\al)}+N\big\lbrack L_{2\al+4,x}^{\al,N}f(x)\big\rbrack\big\vert_{x=0}\,g(0)=\frac{1}{\al!} \int_{0}^{\infty}D_x \big\lbrack e^{-x}x^{\al +1} D_xf(x)\big\rbrack g(x)dx\\  
  &+\frac{N(-1)^{\al+1}}{(\al+2)!\,\al!}\int_{0}^{\infty}D_x^{\al +2} \big\lbrace
  e^{-x} D_x^{\al +2}\big\lbrack x^{\al +1}f(x)\big\rbrack \big\rbrace
  x^{\al+1}g(x)dx +N(\al+1)f'(0)g(0). 
  \end{aligned}
  \label{eq4.3}	 
  \end{equation}  
Integrating the first integral on the right-hand side by parts, the integrated term clearly vanishes at the origin and at infinity, and we are left with 
 \begin{equation*}
  -\frac{1}{\al!} \int_{0}^{\infty}e^{-x}x^{\al +1}f'(x)g'(x)dx
   \end{equation*}  
As to the second integral, an $(\al+2)$-fold integration by parts yields
 \begin{equation*}
  \begin{aligned}     
  \frac{N}{(\al+2)!\,\al!}\,\bigg\lbrace &
  \sum_{j=0}^{\al+1}(-1)^{j+\al+1} D_x^{\al +1-j} \big\lbrace
    e^{-x} D_x^{\al+2}\big\lbrack x^{\al +1}f(x)\big\rbrack \big\rbrace
  D_x^j\big\lbrack x^{\al +1}g(x)\big\rbrack \vert_{x=0}^\infty\\
  &-\int_{0}^{\infty}e^{-x} D_x^{\al+2}\big\lbrack x^{\al +1}f(x)\big\rbrack 
   D_x^{\al+2}\big\lbrack x^{\al +1}g(x)\big\rbrack dx\bigg\rbrace.  
   \end{aligned}
   \end{equation*}  
Since all terms of the sum vanish up to the one for $j=\al+1$, evaluated at $x=0$, its contribution is just
\begin{equation*}
 -\frac{N}{(\al+2)!\,\al!} D_x^{\al+2}\big\lbrack x^{\al +1}f(x)\big\rbrack 
 D_x^{\al+1}\big\lbrack x^{\al +1}g(x)\big\rbrack \vert_{x=0}=
  -N(\al +1)f'(0)g(0).
    \end{equation*} 
Putting all parts together, the two constant terms compensate, while the remaining intgral expressions are symmetric with respect to the functions $f,g$. Hence, the scalar product (4.3) coincides with the right-hand side of identity  (\ref{eq4.2}).  
\hfill $\square$
\end{Proof}
\begin{corollary}
\label{cor4.2}
Let $\al \in \mathbb{N}_0$ and $N > 0$. Given the Laguerre-type differential equation  (\ref{eq4.1}), the following two features are equivalent.\\
(i) The differential operator $L_{2\al+4,x}^{\al,N}$ satisfies the symmetry relation (\ref{eq4.2}). \\
(ii) The Laguerre-type polynomials $y_n(x)=L_n^{\al,N}(x)$, $n\in\mathbb{N}_0$, fulfill the orthogonality relation 
\begin{equation}
\big(y_n,y_m\big)_{w(\al,N)}=h_n\,\delta _{n,m},\; n,m \in \mathbb{N}_0.
  \label{eq4.4} 
  \end{equation}
\end{corollary}
\begin{Proof}
In view of Eq.(\ref{eq4.1}), the symmetry property (\ref{eq4.2}) of the operator $L_{2\al+4,x}^{\al,N}$ implies that
\begin{equation}
  \begin{aligned}
 &\big(\La_{2\al+4,n}^{\al,N}-\La_{2\al+4,m}^{\al,N}\big)     
 \,(y_n,y_m)_{w(\al,N)}\\
 &=-\big(L_{2\al+4,x}^{\al,N}y_n,y_m\big)_{w(\al,N)}+
 \big(y_n,L_{2\al+4,x}^{\al,N}y_m\big)_{w(\al,N)}=0.
 \end{aligned}  
 \label{eq4.5} 
 \end{equation}
Since the eigenvalues $\La_{2\al+4,n}^{\al,N}$ are strictly increasing for $n \in \mathbb{N}_0$, the scalar product (\ref{eq4.4}) vanishes for all $n \ne m$. 

Vice versa, a combination of Eq.(\ref{eq4.1}) and the orthogonality relation (\ref{eq4.4}) readily yields, for any two polynomial $p_n(x), q_m(x)$ of degree $n$ and $m$, respectively, that  
\begin{equation}
  \begin{aligned}
\big(L_{2\al+4,x}^{\al,N}p_n,q_m\big)_{w(\al,N)}&=
\sum_{k=0}^{n}\sum_{j=0}^{m}(p_n,y_k)_w(q_m,y_j)_w\, h_k^{-1} h_j^{-1}\big(L_{2\al+4,x}^{\al,N}y_k,y_j\big)_w\\
&=-\sum_{k=0}^{\min{(n,m)}}(p_n,y_k)_w(q_m,y_k)_w \, h_k^{-1} \La_{2\al+4,k}^{\al,N}\\
&=\big(p_n,L_{2\al+4,x}^{\al,N}q_m\big)_{w(\al,N)}
  \end{aligned}
  \end{equation}
Since the differential operator is bounded and the polynomials are dense in the corresponding function space, the symmetry relation (4.2) holds in general.
\hfill $\square$
    \end{Proof}

 \section{Elementary representation of the Bessel-type equation}
 \label{sec:5}
     
In 1994, W. N. Everitt and the author \cite{EM} introduced the continuous system of Bessel-type functions $\{J_{\la}^{\al,M}(x)\} _{\la \ge 0}$, $\al \in \mathbb{N}_0$, $M \ge 0$, $0 \le x< \infty$, which satisfy an orthogonality relation in distributional sense and thus give rise to a generalized Hankel transform and its inverse \cite{EKLM}. Defined via the following limit of the Laguerre-type polynomials,
\begin{equation}
 J_{\la}^{\al,M}(x)=\lim_{n \rightarrow \infty}\frac{n!}{(\al +1)_n} L_n^{\al ,N(n)}
 \bigg (\frac{(\la x)^2}{4n}\bigg),\;N(n)=\frac{(\la/2)^{2\al+2}}{(n+1)_{\al+1}}M,
      \label{eq5.1}   
      \end{equation}
the Bessel-type functions are explicitly given by
\begin{equation}
 J_\la^{\al,M}(x)=J_\la^\al(x)+M\,K_\la^\al(x),\;
 K_\la^\al(x)=- k_\la^\al x^2 J_\la^{\al+2}(x),\; 
  k_\la^\al =\frac{(\la/2)^{2\al+4}}{(\al+1)(\al+2)!}.
      \label{eq5.2}   
      \end{equation}
Here, 
\begin{equation*}
 J_\la^\al(x)={}_0 F_1(-;\al+1;-\tfrac{1}{4}(\la x)^2)=
 2^\al\Ga(\al+1)(\la x)^{-\al} J_\al(\la x),\; \la \ge 0,
 \end{equation*}	
denote the classical Bessel functions associated with the differential equation
\begin{equation}
\begin{aligned}
\big\lbrack D_x^2 +\frac{2\al+1}{x}D_x\big\rbrack J_\la^\al(x)= -\la^2 J_\la^\al(x).
\end{aligned}
 \label{eq5.3}   
 \end{equation}
 
\begin{theorem}
\label{thm5.1}
(i) For $\al \in \mathbb{N}_0,\, M \ge 0$, the Bessel-type functions $J_\la^{\al,M}(x), \la \ge 0$, satify the eigenvalue equation
\begin{equation}
\bigg\lbrace \big\lbrack \widetilde{L}_{2,x}^\al +\la^2\big\rbrack +
\frac{M}{2^{2\al+2}(\al+2)!}\big\lbrack \widetilde{L}_{2\al+4,x}^\al+\la^{2\al+4}
\big\rbrack \bigg\rbrace J_\la^{\al,M}(x)=0,\;0 \le x<\infty,
 \label{eq5.4}   
 \end{equation}
where, for sufficiently smooth functions $y(x)$,
\begin{equation}
\begin{aligned}
\widetilde{L}_{2,x}^\al y(x)=&x^{-2\al-1}D_x \big\lbrack x^{2\al+1}D_x y(x)\big\rbrack=\big\lbrack D_x^2 +\frac{2\al+1}{x}D_x\big\rbrack y(x),\\
\widetilde{L}_{2\al+4,x}^\al y(x)=&(-1)^{\al +1}x^2(x^{-1}D_x)^{2\al+4} \lbrack x^{2\al+2} y(x)\rbrack\\
=&(-1)^{\al +1}\bigg\lbrack D_x^2 +\frac{2\al+1}{x}D_x-\frac{4\al +4}{x^2} \bigg\rbrack ^{\al +2} y(x).
\end{aligned}
 \label{eq5.5}   
 \end{equation}
(ii) For all $\la \ge 0$, the second component $K_\la^\al(x)$ of the Bessel-type function (5.2) is an eigenfunction of the equation
\begin{equation}
\big\lbrack \widetilde{L}_{2\al+4,x}^\al+\la^{2\al+4}\big\rbrack  K_\la^\al(x)=0,\;0 \le x<\infty.
 \label{eq5.6}   
 \end{equation}   
\end{theorem} 
\begin{Proof}
(i) Consider the Laguerre-type equation (4.1) with the new independent variable $\xi$ and  $N=N(n)$ as chosen in item (5.1). After multiplication by $\la ^2/n$,  the functions $y_n(\xi)=$ $[n!/(\al+1)_n]\,L_n^{\al,N(n)}(\xi),\;0 \le \xi <\infty$, $n \in \mathbb{N}_0$, satisfy
\begin{equation*}
\bigg\lbrace \bigg\lbrack \frac{\la^2}{n}L_{2,\xi}^\al+\la^2\bigg\rbrack +
\frac{M}{2^{2\al+2}(\al+2)!}\bigg\lbrack \frac{\la^{2\al+4}}{(n)_{\al+2}}
L_{2\al+4,\xi}^\al+\la^{2\al+4}\bigg\rbrack \bigg\rbrace y_n(\xi)=0.
 \end{equation*}
Now we substitute $\xi=(\la x)^2/(4n)$ and replace $D_\xi$ by $(2n/\la^2)\,\de_x, \de_x=x^{-1}D_x$. In the limit $n \to \infty$, we readily obtain that
\begin{equation*}
\begin{aligned}
\frac{\la^2}{n}L_{2,\xi}^\al\,y_n(\xi)&=\frac{\la^2}{n}\bigg\lbrack 
\xi\,D_\xi^2 +(\al+1-\xi)D_\xi\bigg\rbrack y_n(\xi)\\
&=\bigg\lbrack x^2\de_x^2+\bigg(\al+1-\frac{(\la x)^2}{4n}\bigg)2\de_x \bigg\rbrack y_n\bigg(\frac{(\la x)^2}{4n}\bigg)\\
&\xrightarrow[n \to \infty]{} \bigg\lbrack D_x^2+\frac{2\al+1}{x}D_x \bigg\rbrack 
J_\la^{\al,M}(x)=\widetilde{L}_{2,x}^\al J_\la^{\al,M}(x).
  \end{aligned}
 \end{equation*}
Furthermore, since $\lim_{n \to \infty}\exp(\pm \xi)=\lim_{n \to \infty}\exp(\pm (\la x)^2/(4n))=1$,
\begin{equation*}
\begin{aligned}
\frac{\la^{2\al+4}}{(n)_{\al+2}}&L_{2\al+4,\xi}^\al y_n(\xi)
=\frac{\la^{2\al+4}}{(n)_{\al+2}}(-1)^{\al+1}e^\xi \xi\,D_\xi^{\al+2}\bigg\lbrace 
e^{-\xi}D_\xi^{\al+2}\big\lbrack \xi^{\al+1} y_n(\xi)\big\rbrack\bigg\rbrace\\
&=\frac{\la^{2\al+4}}{(n)_{\al+2}} \frac{(\la x)^2}{4n}
\bigg(\frac{2n}{\la^2}\bigg)^{2\al+4}\bigg(-\frac{\la^2}{4n}\bigg)^{\al+1}
e^\xi \,\de_x^{\al+2}\bigg\lbrace e^{-\xi}\de_x^{\al+2}\bigg\lbrack x^{2\al+2}
y_n \bigg(\frac{(\la x)^2}{4n}\bigg)\bigg\rbrack\bigg\rbrace\\
&\xrightarrow[n \to \infty]{} (-1)^{\al+1} x^2\de_x^{2\al+4}\big\lbrack x^{2\al+2} 
J_\la^{\al,M}(x)\big\rbrack.
  \end{aligned}
 \end{equation*}
This yields the first representation of $\widetilde{L}_{2\al+4,x}^\al  J_\la^{\al,M}(x)$ in item (5.5). To prove the second line, we first invoke the relation given in \cite [(3.5)]{EM}, 
\begin{equation}
\de_x^{2\be}\lbrack x^{2\be}y(x)\rbrack=
\big\lbrack D_x^2+\frac{2\be+1}{x}D_x \big\rbrack^{\be} y(x),\,y\in C^{(2\be)}(0,\infty), \,\be \in\mathbb{N}_0,
 \label{eq5.7}
\end{equation}
with $\be=\al+2$ and $y(x)=x^{-2}J_{\la}^{\al,M}(x)$. The assertion then follows by applying $\al+2$ times the obvious identity
\begin{equation}
x^2 \bigg\lbrack D_x^2 +\frac{2\al+5}{x}D_x\bigg\rbrack \lbrack x^{-2} y(x)\rbrack=\bigg\lbrack D_x^2 +\frac{2\al+1}{x}D_x-\frac{4\al +4}{x^2} \bigg\rbrack y(x). 
  \label{eq5.8}
  \end{equation} 
	 
(ii) Using the identity (5.8) once more and employing the Bessel equation (5.3) with parameter $\al+2$, we find that
\begin{equation*}
\begin{aligned}
\widetilde{L}_{2\al+4,x}^\al K_\la^{\al}(x)=&
(-1)^{\al +2}k_{\la}^{\al}\bigg\lbrack D_x^2 +\frac{2\al+1}{x}D_x-\frac{4\al +4}{x^2} \bigg\rbrack ^{\al +2} \lbrack x^2 J_{\la}^{\al +2}(x)\rbrack\\
=&(-1)^{\al +2}k_{\la}^{\al}x^2\bigg\lbrack D_x^2 +\frac{2\al+5}{x}D_x\bigg\rbrack ^{\al +2} J_{\la}^{\al +2}(x)\\
=&k_{\la}^{\al} x^2\la^{2\al+4}J_{\la}^{\al+2}(x)=-\la^{2\al+4}K_{\la}^{\al}(x).
 \end{aligned}
  \end{equation*}
Alternatively, we can utilize the differential formulas for the Bessel functions
\begin{equation}
\begin{aligned}
\de_x J_\la^\ga(x)&=-\frac{\la^2}{2(\ga+1)} J_\la^{\ga+1}(x),\;\ga >-1,\\
\de_x\lbrack x^{2\ga} J_\la^\ga(x)\rbrack&=2\ga\,x^{2\ga-2}
J_\la^{\ga-1}(x),\;\ga >0
 \end{aligned}
  \label{eq5.9}
  \end{equation}
to get
\begin{equation*}
\begin{aligned}
\widetilde{L}_{2\al+4,x}^\al K_\la^\al(x)&=(-1)^{\al+2}k_\la^\al x^2\delta_x^{2\al+4}
\big\lbrack x^{2\al+4} J_\la^{\al+2}(x)\big\rbrack\\
&=(-1)^{\al+2}k_\la^\al x^2 2^{\al+2}(\al+2)!\,\de_x^{\al+2} J_\la^0(x)\\
&=k_\la^\al x^2 \la^{2\al+4}J_\la^{\al+2}(x)=-\la^{2\al+4}K_\la^\al(x).\hspace{2cm} \square 
    \end{aligned}
   \end{equation*}    
\end{Proof}  	 
\\
\textbf{Remark.} If we proceed from the first factorization of the Laguerre-type operator in (3.3) and apply the limiting process used in the proof of Theorem 5.1 to each of the $\al +2$   factors, we arrive again at the second identity of the operator $\widetilde{L}_{2\al+4,x}^\al$ in (5.5).\\
 
In \cite[Thm. 2.4]{EM}, the higher-order differential operator in Eq. (5.4) has been given in the form
\begin{equation}
\begin{aligned}
&\widetilde{L}_{2\al+4,x}^\al y(x)=(-1)^{\al+1}\sum_{i=1}^{2\al+4}A_i^\al x^{i-2\al -4} D_x^i y(x),\\
&A_i^\al =\frac{(\al+1)!}{(i-1)!}\sum_{j=\max(i,\al+3)}^{2\al+4}(-1)^{i+j}
\binom{2\al+4}{j}\frac{(2j-i-1)!\;2^{i-2j+2\al+4}}{(j-\al-3)!\,(j-i)!}.
    \end{aligned}
  \label{eq5.10}
  \end{equation}    
\begin{corollary}
\label{cor5.2}
The representations (5.5) and (5.10) of the Bessel-type operator are equivalent.
\end{corollary}
\begin{Proof}
Due to the product rule $\de_x[f \cdot g]=\de_x[f] \cdot g = f \cdot \de_x[g]$, it is clear that
\begin{equation*}
\begin{aligned}
(-1)^{\al+1}\widetilde{L}_{2\al+4,x}^\al y(x)&=x^2\de_x^{2\al+4}\big\lbrack
x^{2\al+2}y(x)\big\rbrack \\
&=x^2\sum_{j=0}^{2\al+4}\binom{2\al+4}{j}\de_x^{2\al+4-j}\big\lbrack x^{2\al+2}
\big\rbrack \de_x^j y(x)\\
&=\sum_{j=\al+3}^{2\al+4}\binom{2\al+4}{j}\frac{2^{2\al+4-j}(\al+1)!}{(j-\al-3)!}x^{2j-2\al-4}\de_x^j y(x).
    \end{aligned}
   \end{equation*}   
Inserting the expansion formula \cite[(2.8)]{EM}, i.e.
\begin{equation*}
\begin{aligned}
x^{2j-2\al-4}\de_x^j y(x)=\sum_{i=1}^{j}\frac{(-2)^{i-j}(2j-i-1)!}
{(j-i)!\,(i-1)!} x^{i-2\al-4}D_x^i y(x),
    \end{aligned}
   \end{equation*}  
and interchanging the role of the two resulting summations, we end up with the required representation (5.10). 
\hfill $\square$
\end{Proof}

\vskip0.5cm
\begin{footnotesize}
\noindent
C. Markett, Lehrstuhl A f\"ur Mathematik, RWTH Aachen, 52056 Aachen, Germany;
\sPP
E-mail: {\tt markett@matha.rwth-aachen.de}
\end{footnotesize}

\begin{thebibliography}{10}

\bibitem{Ba0} H. Bavinck, 
	{\em  On a linear perturbation of the Laguerre operator}, J. Comput. Appl. Math. 106 (1999), 197--202.
	 
\bibitem{Du} A. J. Dur\'{a}n, 
	{\em  Using $\mathcal{D}$-operators to construct orthogonal polynomials satisfying higher order difference or differential equations}, J. Approx. Theory 174 (2013), 10–-53. 
	
\bibitem{HTF2} A. Erd\'{e}lyi, W. Magnus, F. Oberhettinger and F. G. Tricomi,
  	{\em Higher Transcendental Functions, vol. II}, McGraw-Hill, 1953.
	
\bibitem{EM} W. N. Everitt and C. Markett,
	{\em On a generalization of Bessel 	functions satisfying higher-order differential equations}, J Comput. Appl. Math. 54 (1994), 325--349.

\bibitem{EKLM} W. N. Everitt , H. Kalf, L. L. Littlejohn and C. Markett,
	{\em Forth-order Bessel equation: eigenpackets and a generalized Hankel transfom}, Integral Transforms and Special Functions 17 (2006), 845--862.

\bibitem{ELW} W. N. Everitt, L. L. Littlejohn and R. Wellman,
	{\em The symmetric form of the Koekoek’s Laguerre type differential equation}, J. Comput. Appl. Math. 57 (1995), 115--121.

\bibitem{KK1} J. Koekoek and R. Koekoek, 
	{\em On the differential equation for Koornwinder's generalized Laguerre polynomials}, Proc. Amer. Math. Soc. 112 (1991), 1045--1054.

\bibitem{KK2} J. Koekoek and R. Koekoek, 
	{\em Differential equations for generalized Jacobi polynomials}, J. Comput. Appl. Math. 126 (2000), 1--31.
		
\bibitem{Ko} T. H. Koornwinder, 
	{\em Orthogonal polynomials with weight function $(1-x)^{\alpha}(1+x)^{\beta}+M\delta(x+1)+N\delta(x-1)$}, Canad. Math. Bull. 27(2) (1984), 205--214. 

\bibitem{KLY} H. W. Kwon, L. L. Littlejohn and G. J. Yoon,
	{\em Bochner-Krall orthogonal polynomials}, Proc. International Workshop on Special Functions, Asymptotics, Harmonic Analysis and Mathematical Physics,  World Scientific Publishers, Hong Kong, 1999, 181--193.  

\bibitem{LK} L. L. Littlejohn and A. M. Krall,
	{\em On the classification of differential equations having orthogonal polynomial solutions II}, Ann. Mat. Pura Appl. 4 (1988), 77--102. 
 	
\bibitem{Ma2} C. Markett, 
	{\em An elementary representation of the higher-order Jacobi-type differential equation}, Indagationes Mathematicae (2017), http://dx.doi.org/10.1016/j.indag.2017.06.015. 
 	
\bibitem{Ma3} C. Markett, 
	{\em The higher-order differential operator for the generalized Jacobi polynomials - new representation and symmetry}, arXiv:1704.01764v1 [math.CA] (2017).
	
\bibitem{We} R. Wellman, 
	{\em Self-adjoint representations of a certain sequence of spectral differential equations}, Ph.D. Thesis, Utah State Univ. 5, 1995.
	
\end{thebibliography}
\end{document}